\documentstyle[12pt]{article}
\newcommand{\nat}{{\rm I\!N}}
\begin{document}
\thispagestyle{empty}

\begin{center} {\large \bf A Multiple Integral Explicit Evaluation
Inspired by The Multi-WZ Method} \end{center}

\begin{center} Akalu Tefera\footnote{Department of Mathematics, Temple 
University, Philadelphia, PA 19104.\\ akalu@math.temple.edu}  \end{center}

{\bf Abstract}

We give an identity which is conjectured and proved  by using
an implementation [3] in multi-WZ [5].

\bigskip

{\bf 0. \qquad  Introduction}

There are relatively few known non-trivial evaluations
of n-dimensional integrals,  with {\it arbitrary $n$}.
Celebrated examples are the Selberg and the Metha-Dyson integrals,
as well as the Macdonald constant term ex-conjectures for the various
infinite families of root systems. They are all very important. See [1] for
a superb exposition of the various known proofs and of
numerous intriguing applications.

At present, the (continuous version of the) WZ method[5] is
capable of mechanically proving these identities only
for a fixed $n$. In principle for {\it any} fixed $n$
(even, say, $n=100000$), but in practice only for $n \leq 5$.
However, by interfacing a human to the computer-generated
output, the human may discern a pattern, and generalize the
computer-generated proofs for $n=1,2,3,4$ to an arbitrary
$n$.

Using this strategy, Wilf and Zeilberger[5] gave a WZ-style proof
of Selberg's integral evaluation. But just giving 
yet another proof
of an already known identity, especially one that already had
(at least) three beautiful proofs (Selberg's, Aomoto's, and Anderson's,
see [1]), is not very exciting.

In this article we present a {\it new} multi-integral
evaluation, that was {\it first} found using the author's
implementation of the continuous multi-WZ method[3]. 
Both the conjecturing part, and the proving part, were
done by a close human-machine collaboration. 
Our proof hence may be termed {\it computer-assisted}
but not yet {\it computer-generated}.

Now that the result is known and proved, it may be of interest
to have a non-WZ proof, possibly by performing an appropriate
change of variables, converting the multi-integral to a double
integral. My advisor, Doron Zeilberger, is offering
\$100 for such a proof, provided it does not exceed the 
length of the present proof.

\bigskip

{\bf 1. \qquad  Notation}

\begin{tabular}{ll}
${\bf x} = (x_{1}, \ldots, x_{k})$, & $(y)_{m} = \prod_{i=0}^{m-1}(y + 
i)$,\\ 
$d{\bf x} = dx_{1} \cdots dx_{k}$, & $e_1({\bf x})  = \sum_{i=1}^{k} 
x_{i}$,\\
$\hat {{\bf x}}_{i} = (x_{1}, \ldots, \hat {x}_{i}, \ldots, x_{k})$, & 
$e_2({\bf 
x}) = \sum_{i<j} x_{i} x_{j}$,\\
$\Delta_n F(n, {\bf x}) = F(n + 1, {\bf x}) - F(n, {\bf x})$ &\\
\end{tabular}

\bigskip

{\bf 2. \qquad  The Integral Evaluation}

{\bf Theorem }
\[  \int_{[0, +\infty)^{k}} (e_{2}({\bf 
x}))^{m} (e_{1}({\bf x}))^{n} e^{-e_{1}({\bf x})} d{\bf x} = \frac {m! (2 m  
+ n 
+ k - 1)! {(k/2)}_{m} } {(2 m + k - 1)!} \left(\frac {2 (k - 1)}{k} 
\right)^{m}T_{k}(m)\] for all $k$ in $\nat$, and for all $m$, $n$ in 
${\mathbf{Z}}_{\geq 0}$, where, 
\[T_{k}(m) - T_{k}(m - 1) = \frac {(k (k - 2))^{m} {((k - 1)/2)}_{m}}{(k - 1)^{2 
m}{(k/2)}_{m}} T_{k - 1}(m) \]
for all $k \geq 2$, $T_{1} (m) = 0$, for all $m$ in ${\mathbf{Z}}_{\geq 0}$,  
and $ T_{k} (0) = 1 $ for all $ k \geq 2$.

\bigskip

{\bf 3. \qquad  Proof of the  Integral Evaluation}

If $k = 1 $, then trivially, both sides of the integral equate to zero.
Let $k > 1$ and $A_{k}(m, n)$ be the left side of the integral divided by 
\[ \frac {m! (2 m  + n + k - 1)! {(k/2)}_{m} } {(2 m + k - 1)!} \left(\frac 
{2 
(k - 1)}{k} \right)^{m}.\] We want to show $A_{k}(m, n) = T_{k}(m)$, for all 
$m$, $n$ in ${\mathbf{Z}}_{\geq 0}$. Let 
\[ F_{k}(m, n; {\bf x}) := \frac {(2 m + k - 1)!} {m! (2 m  + n + 
k - 1)! {(k/2)}_{m} }  \left(\frac {k} {2 (k - 1)} \right)^{m} (e_{2}({\bf 
x}))^{m} 
(e_{1}({\bf x}))^{n} e^{-e_{1}({\bf x})} \]
We construct\footnote{for specific $k$, the rational function $R$ is 
obtained 
by using {\bf SMint} [3] and the output is available from {\em 
http://www.math.temple.edu/$\sim$akalu/maplepack/rational1.output}}  
\[ R(u; v_{1}, \ldots, v_{k-1}) := \frac {u}{2 m + n + k},\] with the motive 
that 
\[ {\mbox {(WZ 1)}} \qquad \Delta_{n} F_{k}(m, n; {\bf x}) = - 
\sum_{i=1}^{k} 
 D_{x_{i}}[R(x_{i}; \hat{{\bf x}}_{i}) F_{k}(m, n; {\bf x})].\]
Now, we verify (WZ 1),
\begin{eqnarray*}
\lefteqn{ \frac {F_{k} (m, n + 1; {\bf x}) - F_{k}(m, n; {\bf x}) + 
\sum_{i=1}^{k} 
 D_{x_{i}}[R(x_{i}; \hat{{\bf x}}_{i}) F_{k}(m, n; {\bf x})]
}{F_{k}(m, n; {\bf x})}}\\
 & = & \frac {F_{k} (m, n + 1; {\bf x})}{F_{k}(m, n;{\bf x})} - 1 + 
\sum_{i=1}^{k} 
 D_{x_{i}}[R(x_{i}; \hat{{\bf x}}_{i})] + R(x_{i}; \hat{{\bf 
x}}_{i})D_{x_{i}} 
\left[log \left(F_{k}(m, n; {\bf x})\right) \right]\\ 
& = & \frac {e_{1}({\bf x})}{2 m + n + k} - 1 + \frac {k}{2 m + n + k} + \\
& & \sum_{i=1}^{k} \left(\frac {n}{e_{1}({\bf x})}\frac {x_{i}}{2 m + n + k} 
 + 
\frac{m e_{1}(\hat {{\bf x}}_{i})}{e_{2}({\bf x})} 
\frac {x_{i}}{2 m + n + k} - \frac {x_{i}}{2 m + n + k} \right)\\
& = & \frac {e_{1}({\bf x})}{2 m + n + k} - 1 + \frac {k}{2 m + n + k} + 
\frac 
{n}{2 m + n + k} + \frac{2 m }{2 m + n + k} -\frac {e_{1}({\bf x})}{2 m + n 
+ 
k}\\
& = & 0.
\end{eqnarray*}

Hence, by integrating both sides of (WZ 1) w.r.t $x_{1}, \ldots, x_{k}$ 
over $[0, \infty)^{k}$, we get \[ A_{k}(m, n + 1) - A_{k}(m, n) \equiv 0.\]

To complete the proof we show $A_{k}(m, 0) =  T_{k}(m)$.

To this end, set $A_{k}(m) := A_{k}(m, 0)$ and $F_{k}(m; {\bf x}) := 
F_{k}(m, 
0; {\bf x})$. Now, we construct\footnote{
for specific $k$, the rational function $R$ is obtained by using {\bf SMint} 
[3] 
and the output is available from {\em 
http://www.math.temple.edu/$\sim$akalu/maplepack/rational2.output}},
\[ R(u; v_{1}, \ldots, v_{k-1}) := \frac {((k - 1) (m + 1) + e_{1}(v_{1}, 
\ldots, v_{k-1})) u + e_{2}(v_{1}, \ldots, v_{k-1})}{(k - 1)(m + 1)(2 m + 
k)}\] 
with the motive that 
\[ {\mbox {(WZ 2)}} \qquad  F_{k}(m + 1; {\bf x}) - F_{k}(m; {\bf x}) = - 
\sum_{i=1}^{k}  D_{x_{i}}[R(x_{i}; \hat{{\bf x}}_{i}) F_{k}(m; {\bf x})].\]
Verification of (WZ 2):
\begin{eqnarray*}
\lefteqn{ \frac {F_{k}(m + 1; {\bf x}) - F_{k}(m; {\bf x}) + \sum_{i=1}^{k}  
D_{x_{i}}[R(x_{i}; \hat{{\bf x}}_{i}) F_{k}(m; {\bf x})]
}{F_{k}(m; {\bf x})}}\\
 & = & \frac {F_{k} (m + 1; {\bf x})}{F_{k}(m; {\bf x})} - 1 + 
\sum_{i=1}^{k} 
D_{x_{i}}[R(x_{i}; \hat{{\bf x}}_{i})] + \sum_{i=1}^{k} 
R(x_{i}; \hat{{\bf x}}_{i})D_{x_{i}} \left[log \left(F_{k}(m; {\bf 
x})\right) 
\right]\\ 
& = & \frac{k e_{2}({\bf x})}{(m + 1) (k - 1) (2 m + k)} - 1 + 
\sum_{i=1}^{k} 
\frac {(k - 1) (m + 1) + e_{1}(\hat{{\bf x}}_{i})}{(m + 1) (k - 1) (2 m + 
k)} + 
\\
& &\sum_{i=1}^{k} \frac {(k - 1) (m + 1) x_{i} + e_{2}({\bf x}))}{(m + 1) (k 
- 
1) (2 m + k)} \left (\frac {m e_{1}({\hat {\bf 
x}}_{i})}{e_{2}({\bf x})} - 1 \right)\\
& =  & \frac{k e_{2}({\bf x})}{(m + 1) (k - 1) (2 m + k)} - 1 + \frac{k}{2 m 
+ 
k} + \frac {e_{1}({\bf x})}{(m + 1)(2 m + k)} + \frac {2 m}{2 m + k} - \frac 
{e_{1}({\bf x})}{2 m + k} + \\ 
& & \frac {m e_{1}({\bf x})}{(m + 1)(2 m + k)} -\frac {k e_{2}({\bf x})}{(m 
+ 
1)(k - 1)(2 m + k)}\\
& = & 0.
\end{eqnarray*}

Hence, by integrating both sides of (WZ 2) w.r.t. $x_{1}, \ldots, x_{k}$ 
over $[0, \infty)^{k}$, we obtain,

\[A_{k}(m + 1) - A_{k}(m) = \frac {(k (k - 2))^{m + 1} {((k - 1)/2)}_{m + 
1}}{(k - 1)^{2 (m + 1)}{(k/2)}_{m + 1}} A_{k - 1}(m + 1),\]

and noting that $A_{k}(0) = 1$, $A_{1} (m) = 0$, it follows that $A_{k}(m) = 
T_{k} (m)$, for all $m$ in ${\mathbf{Z}}_{\geq 0}$.  Consequently, $A_{k}(m, 
n) = T_{k} (m)$ for all $m$, $n$ in ${\mathbf{Z}}_{\geq 0}$. \hspace*{40mm} 
$ \Box$

By unfolding  the recurrence equation for $T_{k}(m)$, we obtain the following 
identity.

{\bf Corollary}
\begin{eqnarray*}
\lefteqn{  \int_{[0, +\infty)^{k}} (e_{2}({\bf 
x}))^{m} (e_{1}({\bf x}))^{n} e^{-e_{1}({\bf x})} d{\bf x} =\frac {m! (2 m  
+ n + k - 1)! {(k/2)}_{m} } {(2 m + k - 1)!} \left(\frac {2 (k - 1)}{k} 
\right)^{m}}\\
& & \left( 1 + \sum_{r = 1}^{k - 2}\sum_{1 \leq s_{r} \leq \cdots \leq s_{1} 
\leq m}  \prod_{i = 1}^r \frac {((k - i)^2 - 1)^{s_{i}} ((k 
- i)/2)_{s_{i}}}{(k - i)^{2 s_{i}} ((k - i + 1)/2)_{s_{i}}} \right) 
\end{eqnarray*}

\newpage
{\bf 4. \qquad Remarks} 

{\bf 1.} 
\qquad From the computational point of view, the recurrence form of the integral 
is {\em nicer} than  its indefinite summation form (the above corollary),
for the former requires  $O(m k)$ calculations, whereas 
the latter requires $O(m^{k})$ calculations. However, in both forms the 
evaluation of the right side of the integral is much faster (for
specific $m$, $n$, and $k$)
than the direct evalution of the left side of our intergal. 
Hence both forms are indeed complete {\em answers} in the sense of Wilf[4].

{\bf 2.} \qquad The present paper is an example of what Doron Zeilberger[6] 
calls
{\it WZ Theory, Chapter 1 1/2}. Even though, at present, our
proof, for general $n$, was human-generated, it looks almost
computer-generated. It seems that by using John Stembridge's[2]
Maple package for symmetric functions, SF, or an extension of it,
it should be possible to write a new version of {\tt SMint} that
should work for {\it symbolic}, i.e. arbitrary, $n$, thereby
fulfilling the hope raised in [6].

\bigskip

{\bf Acknowledgement:} I thank Doron Zeilberger, my Ph.D. thesis advisor, 
for very helpful suggestions and valuable support.  

\bigskip

{\bf References}

[1] \qquad G. Andrews, R. Askey, and R. Roy, {
\it Special Functions}, Cambridge University \\ \hspace*{13mm} Press, 1998.

[2] \qquad J.R. Stembridge, {\it A Maple package for symmetric functions}, J. \\ 
\hspace*{13mm} Symbolic Comput., {\bf 20}(1995), 755-768.

[3] \qquad A. Tefera, {\bf SMint}{\it (A Maple Package for Multiple 
Integrals)},\\ \hspace*{13mm} {\tt 
http://www.math.temple.edu/$\sim$akalu/maplepack/SMint}.

[4]  \qquad H.S. Wilf,{\it What is an answer?}, Amer. Math. Monthly, {\bf 89} 
(1982), 289-292.

[5] \qquad  H.S. Wilf and D. Zeilberger, {\it An Algorithmic proof theory 
for hypergeometric \\ \hspace*{13mm} (ordinary and "q") multisum/integral 
identities}, Invent. Math., {\bf 108} (1992), \\\hspace*{13mm} 575-633.

[6] \qquad D. Zeilberger, {\it WZ Theory, Chapter II},
The Personal Journal of S.B. Ekhad \\ \hspace*{13mm} and D. Zeilberger,
{\tt http://www.math.temple.edu/$\sim$zeilberg/pj.html}.

\end{document}